\documentstyle[12pt]{article}
\catcode`@=11

\newskip\plaincentering \plaincentering=0pt plus 1000pt minus 1000pt
\def\@plainlign{\tabskip=0pt\everycr={}}
\def\eqalignno#1{\displ@y \tabskip\plaincentering
  \halign to\displaywidth{\hfil$\@lign\displaystyle{##}$\tabskip\z@skip
    &$\@lign\displaystyle{{}##}$\hfil\tabskip\plaincentering
    &\llap{$\@lign##$}\tabskip\z@skip\crcr
    #1\crcr}}
\def\leqalignno#1{\displ@y \tabskip\plaincentering
  \halign to\displaywidth{\hfil$\@lign\displaystyle{##}$\tabskip\z@skip
    &$\@lign\displaystyle{{}##}$\hfil\tabskip\plaincentering
    &\kern-\displaywidth\rlap{$\@lign##$}\tabskip\displaywidth\crcr
    #1\crcr}}
\def\plainLet@{\relax\iffalse{\fi\let\\=\cr\iffalse}\fi}
\def\plainvspace@{\def\vspace##1{\noalign{\vskip##1}}}

\def\intic@{\mathchoice{\hskip5\p@}{\hskip4\p@}{\hskip4\p@}{\hskip4\p@}}
\def\negintic@
 {\mathchoice{\hskip-5\p@}{\hskip-4\p@}{\hskip-4\p@}{\hskip-4\p@}}
\def\intkern@{\mathchoice{\!\!\!}{\!\!}{\!\!}{\!\!}}
\def\intdots@{\mathchoice{\cdots}{{\cdotp}\mkern1.5mu
    {\cdotp}\mkern1.5mu{\cdotp}}{{\cdotp}\mkern1mu{\cdotp}\mkern1mu
      {\cdotp}}{{\cdotp}\mkern1mu{\cdotp}\mkern1mu{\cdotp}}}
\newcount\intno@
\def\iint{\intno@=\tw@\futurelet\next\ints@}
\def\iiint{\intno@=\thr@@\futurelet\next\ints@}
\def\iiiint{\intno@=4 \futurelet\next\ints@}
\def\idotsint{\intno@=\z@\futurelet\next\ints@}
\def\ints@{\findlimits@\ints@@}
\newif\iflimtoken@
\newif\iflimits@
\def\findlimits@{\limtoken@false\limits@false\ifx\next\limits
 \limtoken@true\limits@true\else\ifx\next\nolimits\limtoken@true\limits@false
    \fi\fi}
\def\multintlimits@{\intop\ifnum\intno@=\z@\intdots@
  \else\intkern@\fi
    \ifnum\intno@>\tw@\intop\intkern@\fi
     \ifnum\intno@>\thr@@\intop\intkern@\fi\intop}
\def\multint@{\int\ifnum\intno@=\z@\intdots@\else\intkern@\fi
   \ifnum\intno@>\tw@\int\intkern@\fi
    \ifnum\intno@>\thr@@\int\intkern@\fi\int}
\def\ints@@{\iflimtoken@\def\ints@@@{\iflimits@
   \negintic@\mathop{\intic@\multintlimits@}\limits\else
    \multint@\nolimits\fi\eat@}\else
     \def\ints@@@{\multint@\nolimits}\fi\ints@@@}
\def\Sb{_\bgroup\vspace@
        \baselineskip=\fontdimen10 \scriptfont\tw@
        \advance\baselineskip by \fontdimen12 \scriptfont\tw@
        \lineskip=\thr@@\fontdimen8 \scriptfont\thr@@
        \lineskiplimit=\thr@@\fontdimen8 \scriptfont\thr@@
        \Let@\vbox\bgroup\halign\bgroup \hfil$\scriptstyle
            {##}$\hfil\cr}
\def\endSb{\crcr\egroup\egroup\egroup}
\def\Sp{^\bgroup\vspace@
        \baselineskip=\fontdimen10 \scriptfont\tw@
        \advance\baselineskip by \fontdimen12 \scriptfont\tw@
        \lineskip=\thr@@\fontdimen8 \scriptfont\thr@@
        \lineskiplimit=\thr@@\fontdimen8 \scriptfont\thr@@
        \Let@\vbox\bgroup\halign\bgroup \hfil$\scriptstyle
            {##}$\hfil\cr}
\def\endSp{\crcr\egroup\egroup\egroup}
\def\Let@{\relax\iffalse{\fi\let\\=\cr\iffalse}\fi}
\def\vspace@{\def\vspace##1{\noalign{\vskip##1 }}}
\def\aligned{\,\vcenter\bgroup\plainvspace@\plainLet@\openup\jot\m@th\ialign
  \bgroup \strut\hfil$\displaystyle{##}$&$\displaystyle{{}##}$\hfil\crcr}
\def\endaligned{\crcr\egroup\egroup}
\def\matrix{\,\vcenter\bgroup\plainLet@\plainvspace@
    \normalbaselines
  \m@th\ialign\bgroup\hfil$##$\hfil&&\quad\hfil$##$\hfil\crcr
    \mathstrut\crcr\noalign{\kern-\baselineskip}}
\def\endmatrix{\crcr\mathstrut\crcr\noalign{\kern-\baselineskip}\egroup
                \egroup\,}
\newtoks\hashtoks@
\hashtoks@={#}
\def\format{\crcr\egroup\iffalse{\fi\ifnum`}=0 \fi\format@}
\def\format@#1\\{\def\preamble@{#1}%
  \def\c{\hfil$\the\hashtoks@$\hfil}%
  \def\r{\hfil$\the\hashtoks@$}%
  \def\l{$\the\hashtoks@$\hfil}%
  \setbox\z@=\hbox{\xdef\Preamble@{\preamble@}}\ifnum`{=0 \fi\iffalse}\fi
   \ialign\bgroup\span\Preamble@\crcr}

\def\cases{\left\{\,\vcenter\bgroup\plainvspace@
     \normalbaselines\openup\jot\m@th
      \plainLet@\ialign\bgroup$\displaystyle{##}$\hfil&\quad$\displaystyle{{}##}$\hfil\crcr
      \mathstrut\crcr\noalign{\kern-\baselineskip}}

\newif\iftagsleft@
\tagsleft@true
\def\TagsOnRight{\global\tagsleft@false}
\def\tag#1$${\iftagsleft@\leqno\else\eqno\fi
 \hbox{\def\pagebreak{\global\postdisplaypenalty-\@M}%
 \def\nopagebreak{\global\postdisplaypenalty\@M}\rm(#1\unskip)}%
  $$\postdisplaypenalty\z@\ignorespaces}
\interdisplaylinepenalty=\@M
\def\plainallowdisplaybreak@{\def\allowdisplaybreak{\noalign{\allowbreak}}}
\def\plaindisplaybreak@{\def\displaybreak{\noalign{\break}}}
\def\align#1\endalign{\def\tag{&}\plainvspace@\plainallowdisplaybreak@\plaindisplaybreak@
  \iftagsleft@\plainlalign@#1\endalign\else
   \plainralign@#1\endalign\fi}
\def\plainralign@#1\endalign{\displ@y\plainLet@\tabskip\plaincentering\halign to\displaywidth
     {\hfil$\displaystyle{##}$\tabskip=\z@&$\displaystyle{{}##}$\hfil
       \tabskip=\plaincentering&\llap{\hbox{\rm(##\unskip)}}\tabskip\z@\crcr
             #1\crcr}}
\def\plainlalign@
 #1\endalign{\displ@y\plainLet@\tabskip\plaincentering\halign to \displaywidth
   {\hfil$\displaystyle{##}$\tabskip=\z@&$\displaystyle{{}##}$\hfil
   \tabskip=\plaincentering&\kern-\displaywidth
        \rlap{\hbox{\rm(##\unskip)}}\tabskip=\displaywidth\crcr
               #1\crcr}}

\def\re@#1{\par\hangindent\parindent\indent\llap{#1\enspace}\ignorespaces}
\def\qfootnote#1{\edef\@sf{\spacefactor\the\spacefactor}{}#1\@sf
      \insert\footins{\let\egroup=}\footnotesize % latex
      \interlinepenalty100 \let\par=\endgraf
        \leftskip=0pt \rightskip=0pt
        \splittopskip=10pt plus 1pt minus 1pt \floatingpenalty=20000
   \smallskip\re@{#1}\bgroup\strut\aftergroup{\strut\egroup}\let\next}
\catcode`\@=\active
\TagsOnRight
\begin {document}
\makeatletter

\begin{center}
{\large\bf Ultimate Generalization to Monotonicity for Uniform

Convergence of Trigonometric Series}

\vspace{5mm}

Song Ping Zhou ,\footnote{Corresponding author. W. F. James Chair
Professor of St. Francis Xavier University. Research also
supported in part by NSF of China under grant number 10471130.}\\
{\it Institute of Mathematics, Zhejiang  Sci-Tech University,
Xiasha Economic Development Area,
Hangzhou, Zhejiang 310018  China, and\\
Department of Mathematics, Statistics \& Computer Science,
St. Francis Xavier University, Antigonish, Nova Scotia, Canada B2G 2W5\\
szhou@zjip.com }

Ping Zhou,\footnote{Research suppoted by NSERC of Canada} Dan
Sheng Yu \footnote{Research supported in part by NSERC RCD grant
of St. Francis Xavier University and in part by AARMS of Canada}\\
{\it Department of Mathematics, Statistics \& Computer Science,
St. Francis Xavier University, Antigonish, Nova Scotia, Canada B2G 2W5\\
pzhou@stfx.ca,  dyu@stfx.ca}

\end{center}

{\bf Abstract}\\

 Chaundy and Jolliffe [4] proved that if
$\{a_{n}\}$ is a non-increasing (monotonic) real sequence with
$\lim\limits_{n\rightarrow \infty }a_{n}=0$, then a necessary and
sufficient condition for the uniform convergence of the series
$\sum_{n=1}^{\infty }a_{n}\sin nx$ is $ \lim\limits_{n\rightarrow
\infty }na_{n}=0$. We generalize (or weaken) the monotonic
condition on the coefficient sequence $\{a_{n}\}$ in this
classical result to the so-called mean value bounded variation
condition and prove that the generalized condition cannot be
weakened further. We also establish an analogue to the generalized
Chaundy and Jolliffe theorem in the complex space.\\



\section{Introduction and results}

Let $\{a_{n}\}$ be a nonnegative sequence, write
\begin{equation}
\sum_{n=1}^{\infty }a_{n}\sin nx  \label{a1}
\end{equation}%
as a sine series. Denote by $C_{2\pi }$ the space of all continuous
functions of period $2\pi $ equipped with the norm
$$
\Vert f\Vert =\max_{0\leq x\leq 2\pi }|f(x)|.
$$%
In 1916, Chaundy and Jolliffe [4] proved that if $\{a_{n}\}$ is
a non-increasing (monotonic) real sequence (in symbol, $\{a_{n}\}\in %
\mbox{\rm MS}$, i.e. Monotonic Sequence) with
$\lim\limits_{n\rightarrow \infty }a_{n}=0$, then a necessary and
sufficient condition for the uniform convergence of series (1) is
$\lim\limits_{n\rightarrow \infty }na_{n}=0$. This classical
result, together with other convergence results of series
(\ref{a1}), such as $L^{1}$-convergence, $L^{p}$-convergence, and
best approximation, have had many
applications in analysis and attracted lots of attentions.\\

In order to weaken the monotonic condition on the sequence
$\{a_{n}\}$ for series (1) to be uniformly convergent, several
groups, led by Leindler ([8]-[11]), Stanojevic ([15]-[17]) ),
Telyakovskii ([18]-[21]), S. P. Zhou ([7],[23]-[25]), Belov ([2]),
as well as Tikhonov ([22]), etc., have been working on this
problem in various ways for decades and trying to find the best
way to weaken the monotonic condition of the sequence $\{a_{n}\}$
for the sine
series to be uniformly convergent.\\

One way is to generalize the monotonic condition to the
quasimonotone conditions. The classical definition for a sequence
$\{a_{n}\}$ to be
quasimonotone (in symbol, $\{a_{n}\}\in \mbox{\rm CQMS},$ i.e. \textit{%
Classical Quasi-Montone Sequence}) is that if there is an $\alpha
\geq 0$ such that $a_{n}/n^{\alpha }$ is decreasing for all $n>0$
(see [14],[1],[5],[6],[13] ). The general definition for
quasimonotone is the so-called $O$-regularly varying
quasimonotone condition. Let $R(n)$ be an increasing sequence with $%
R(2n)/R(n)$ bounded for all $n>0$. A sequence $\{a_{n}\}$ is said to be $O$%
\textit{-regularly varying quasimonotone sequence} ($\{a_{n}\}\in
\mbox{\rm RVQMS}$) if for some $R(n)$ with the above properties,
$a_{n}/R(n)$ is
decreasing for all $n>0.$ It is proved that the monotonic condition, $%
\{a_{n}\}\in \mbox{\rm MS}$, in the classical Chaundy-Jollif
Theorem, can be
generalized to $\{a_{n}\}\in \mbox{\rm CQMS}$, or more generally, $%
\{a_{n}\}\in \mbox{\rm RVQMS}$. There are numerous works related
to this
topic, for example, see [3], [12],[15]-[24].\\

Although the RVQMS seems very general, it is almost impossible for
one to prove that a sequence $\left\{ a_{n}\right\} $ (without
missing any terms)
is not a RVQMS. The reason is that one has to prove that the sequence $%
\left\{ a_{n}/R(n)\right\} $ is not decreasing for any $R(n)$ with $%
R(2n)/R(n)$ bounded! This may also be one of the main reasons that
analysts gradually lose their interests towards RVQMS.\\

People then move to another direction to establish a new way of
generalizing the monotonic condition by using the so-called
bounded variation concept. Leindler [8] first raised the rest
bounded variation
condition. A nonnegative sequence $\mbox{\bf A}=\{a_{n}\}$ with $%
\lim\limits_{n\rightarrow \infty }a_{n}=0$ is said to be a \textit{rest
bounded variation sequence} ($\{a_{n}\}\in \mbox{\rm RBVS}$) if
$$
\sum_{k=n}^{\infty }|a_{k}-a_{k+1}|\leq C(\mbox{\bf A})a_{n}
$$
holds for all $n=1,2,\cdots $ and some constant $C(\mbox{\bf A})$
depending only upon the sequence $\mbox{\bf A}$. The
Chaundy-Jollif Theorem is again proved to be true in [8] if we
replace the monotonic condition by the RBV condition. However,
Leindler himself proved that CQMS and RBVS are
not comparable in [9].\\

Very recently, Le and Zhou [7] introduced a condition which
generalizes  both CQMS and RBVS. A nonnegative sequence $\mbox{\bf A}%
=\{a_{n}\}$ is said to be a \textit{group bounded variation sequence} ($%
\{a_{n}\}\in \mbox{\rm GBVS}$) if for some given $N_{0}\geq 1,$
$$
\sum_{k=n}^{2n}|a_{k}-a_{k+1}|\leq C(\mbox{\bf A})\max_{n\leq k\leq
n+N_{0}}a_{k}
$$%
holds for some constant $C(\mbox{\bf A})$ and all $n=1,2,\cdots .$ The
monotonic condition in the Chaundy-Jollif Theorem is then extended to $%
\{a_{n}\}\in \mbox{\rm GBVS}$. Later, Yu and Zhou [25] introduced
further the non-onesided bounded variation condition. A
nonnegative sequence $\mbox{\bf A}=\{a_{n}\}$ is said to be a
\textit{non-onesided bounded variation sequence} ($\{a_{n}\}\in
\mbox{\rm NBVS}$) if
$$
\sum_{k=n}^{2n}|a_{k}-a_{k+1}|\leq C(\mbox{\bf A})(a_{n}+a_{2n})
$$%
holds for some constant $C(\mbox{\bf A})$ and all $n=1,2,\cdots .$ Again the
monotonic condition in the Chaundy-Jollif Theorem is extended further to $%
\{a_{n}\}\in \mbox{\rm NBVS}.$\\

Another recent temptation of generalization is to the so-called almost
monotonic sequence. A nonnegative sequence $\mathbf{b}=\{b_{n}\}$ is said to
be an \textit{almost monotonic sequence} ($\{b_{n}\}\in \mbox{\rm AMS}$) if
there is a positive constant $C(\mbox{\bf b})$ such that
$$
b_{k}\leq C(\mbox{\bf b})b_{n}\;\;\mbox{\rm for all}\;\;k\geq n.
$$
An AMS looks easy to manage. Indeed, AMS contains $\mbox{\rm RVQMS}\cup %
\mbox{\rm RBVS}$, but it is not comparable with GBVS, NBVS or
MVBVS (see[10],[11],[25] for more discussion on this). We prove in
this paper that AMS is not an
option to generalize the Chaundy-Jollif Theorem:\\

%\begin{thm}
{\bf Theorem 1 } \quad {\it There exists a sequence $\{b_{n}\}\in \mbox{\rm AMS}$ with $%
\lim\limits_{n\rightarrow \infty }nb_{n}=0$ such that the series $%
\sum_{n=1}^{\infty }b_{n}\sin nx$ is not uniformly convergent.}\\
%\end{thm}

Our main objective of this paper is to generalize the monotonic condition in
the Chaundy-Jollif Theorem to the so-called \textit{mean value bounded
variation condition} and prove that the generalization achieved in this
paper is final.\\

 %\begin{defn}
  {\bf Definition 2} \quad {\it A nonnegative sequence $\mbox{\bf A}=\{a_{n}\}_{n=0}^{\infty }$ is said to be a \textit{mean value
bounded variation sequence} ($\{a_{n}\}\in \mbox{\rm MVBVS}$) if
there is a $\lambda \geq 2$ such that
$$
\sum_{k=n}^{2n}|a_{k}-a_{k+1}|\leq \frac{C(\mbox{\bf A})}{n}\sum_{k=[\lambda
^{-1}n]}^{[\lambda n]}a_{k}
$$
holds for all $n=1,2,\cdots $ and some constant $C(\mbox{\bf A})$
depending only upon the sequence $\mbox{\bf A}$.}\\
%\end{defn}

From the definition, we can see that a MVBVS can either be
non-increasing almost everywhere, or non-decreasing almost
everywhere, and converge to its limit. We show that the class of
MVBVS contains all known classes of sequences mentioned earlier,
except for the AMS, as following propositions:\\

%\begin{prop}
{\bf Proposition 3} \quad {\it If $\mbox{\bf A}=\{a_{n}\}\in \mbox{\rm GBVS}$ in general sense, i.e., $%
\{a_{n}\}$ satisfies
\begin{equation}
\sum_{k=n}^{2n}|a_{k}-a_{k+1}|\leq C(\mbox{\bf A})\max_{n\leq
k<n+N_{0}}a_{k}  \label{a3}
\end{equation}%
for some given $N_{0}\geq 1$, then $\{a_{n}\}\in \mbox{\rm
MVBVS}$. But the reverse is not true, i.e. there are sequences in
$\mbox{\rm MVBVS}$ not satisfying (2).}\\
%\end{prop}

%\begin{prop}
{\bf Proposition 4} \quad {\it If $\mbox{\bf A}=\{a_{n}\}\in
\mbox{\rm NBVS}$, then $\{a_{n}\}\in \mbox{\rm
MVBVS}$. But the reverse is not true, i.e. there are sequences in $%
\mbox{\rm MVBVS}$ which are not in $\mbox{\rm NBVS}$.}\\
%\end{prop}

Our first main result is that the monotonic condition in the
Chaundy-Jollif Theorem can be weakened to $\{a_{n}\}\in \mbox{\rm
MVBVS}$:\\

%\begin{thm}
{\bf Theorem 5} \quad {\it If $\mbox{\bf A}=\{a_{n}\}\in \mbox{\rm
MVBVS}$, then a necessary and sufficient condition either for the
uniform convergence of series (1), or for the continuity of its sum function $f$, is that $%
\lim\limits_{n\rightarrow \infty }na_{n}=0$.}\\
%\end{thm}

We also prove that the MVBV condition cannot be weakened any
further to guarantee the uniform convergence of the sine series
(1), and therefore $\{a_{n}\}\in \mbox{\rm MVBVS}$ is the ultimate
generalization to the monotonic condition in
Chaundy-Jollif Theorem:\\

%\begin{thm}
{\bf Theorem 6}\quad {\it Let $\{M_{n}\}$ be a given nonnegative
increasing sequence tending to
infinity. Then there exists a sine series of the form (1) with $%
\lim_{n\rightarrow \infty }na_{n}=0$ such that for any given $\lambda \geq 2$
$$
\lim_{n\rightarrow \infty }\frac{\sum_{k=n}^{2n}|\Delta a_{k}|}{\frac{M_{n}}{%
n}\sum\limits_{k=[\lambda ^{-1}n]}^{[\lambda n]}a_{k}}=0,
$$
however, the series is not uniformly convergent. }\\
%\end{thm}

We prove the above propositions and theorems in next section and
establish an analogue of Theorem 5 in the complex space in the
last section. We will also investigate other important classic
results in Fourier analysis under the MVBV condition in separate
papers as continuations to this paper.\\

Finally in this section, we summary the generalization of the monotone
conditions in the following two figures. Figure 1 shows the development of
the generalization successively, while Figure 2 shows the relations of the
different generalized classes of monotonic sequences. Here in Figure 2, for
convenience, GBVS is when $N_{0}=1$ of the general class of GBVS.

\unitlength=0.5mm
\begin{picture}(160,220)\thicklines
\put(95,190){\framebox(30,10){MS}}
\put(110,190){\vector(-3,-1){15.9}} \put(95,185){\line(-3,-1){15}}
\put(65,170){\framebox(30,10){RBVS}}
\put(110,190){\vector(3,-1){15.9}} \put(125,185){\line(3,-1){15}}
\put(125,170){\framebox(30,10){CQMS}}
\put(80,170){\vector(-3,-1){32}} \put(50,160){\line(-3,-1){30}}
\put(5,140){\framebox(30,10){AMS}} \put(20,140){\vector(0,-1){10}}
\put(20,130){\line(0,-1){10}} \put(5,110){\framebox(30,10){STOP}}
\put(140,170){\vector(0,-1){15}} \put(140,152){\circle{1}}
\put(140,149){\circle{1}}\put(140,146){\circle{1}}
\put(140,142){\vector(0,-1){15}} \put(125,117){\framebox(30,10){
RVQMS}}\put(125,122){\vector(-4,1){46.5}}\put(80.4,133.15){\line(-4,1){45}}
\put(140,117){\vector(-1,-1){10}}\put(130,107){\line(-1,-1){12.2}}
\put(104,84.2){\framebox(30,10){GBVS}}
\put(80,170){\vector(1,-2){19}}\put(98.8,132.4){\line(1,-2){18.8}}
\put(119,84.2){\vector(0,-1){10}}\put(119,74.2){\line(0,-1){10}}
\put(104,54.2){\framebox(30,10){NBVS}}
\put(119,54.2){\vector(0,-1){10}}\put(119,44.2){\line(0,-1){10}}
\put(104,24.2){\framebox(30,10){MVBVS}}
\put(80,0){\bf Fig.1}
\end{picture}

\unitlength=1mm
\begin{picture}(120,220)\thicklines
\put(60,110){\circle{40}} \put(57,108){\bf MS}
\qbezier(55,100)(45,115)(55,130) \qbezier(65,100)(75,115)(65,130)
\qbezier(55,130)(60,135)(65,130) \qbezier(55,100)(60,95)(65,100)
\put(54,125){\bf RBVS} \qbezier(55,90)(45,105)(55,120)
\qbezier(65,90)(75,105)(65,120) \qbezier(55,120)(60,125)(65,120)
\qbezier(55,90)(60,85)(65,90) \put(54,93){\bf CQMS}
\qbezier(55,80)(40,110)(55,140) \qbezier(65,80)(80,110)(65,140)
\qbezier(55,140)(60,150)(65,140) \qbezier(55,80)(60,70)(65,80)
\put(54,82){\bf GBVS} \qbezier(50,70)(35,110)(50,150)
\qbezier(70,70)(85,110)(70,150) \qbezier(50,150)(60,170)(70,150)
\qbezier(50,70)(60,50)(70,70) \put(54,68){\bf NBVS}
\qbezier(45,60)(30,110)(45,160) \qbezier(75,60)(90,110)(75,160)
\qbezier(45,160)(60,195)(75,160) \qbezier(45,60)(60,25)(75,60)
\put(52,50){\bf MVBVS} \qbezier(25,95)(0,110)(25,125)
\qbezier(95,95)(120,110)(95,125) \qbezier(25,125)(60,143)(95,125)
\qbezier(25,95)(60,77)(95,95) \put(25,105){\bf AMS}\put(90,105){\bf
AMS} \put(50,35){\bf Fig.2}
\end{picture}

\section{Proofs}

Throughout this paper, we always use $C\left( x\right) $ to denote a
positive constant depending only upon $x$, where $x$ can be numbers or
sequences, and use $C$ to denote an absolute positive constant. $C\left(
x\right) $ or $C$ may have different values in different occurrences.

%\begin{pf}
{\bf PROOF.} [Proof of Theorem 1] Let $n_{1}=1$, $n_{2}=10$,
$n_{j+1}=n_{j}^{2}$ for $j=2,3,\ldots ,$ and let
$$
b_{k}=1,\;\;1\leq k\leq 40.
$$
For $j\geq 2$ and $k=1,2,\ldots ,n_{j}-1,$ let
\begin{eqnarray*}
b_{m} &=&\frac{1}{\sqrt{\log n_{j}}}\frac{1}{m},\;\;4kn_{j}\leq
m<(4k+2)n_{j}, \\
b_{m} &=&\frac{1}{8\sqrt{\log n_{j}}}\frac{1}{m},\;\;(4k+2)n_{j}\leq
m<4(k+1)n_{j}.
\end{eqnarray*}%
Then $nb_{n}\rightarrow 0$, $n\rightarrow \infty $, and $b_{k}\leq 8b_{n}$
for all $k>n$ (this means that $\{b_{n}\}$ is an almost monotonic sequence).
Therefore the series $\sum_{n=1}^{\infty }b_{n}\sin nx$ is well defined.
Denote by $S_{n}(x)$ the $n$th partial sum of the series, i.e.%

%\begin{equation*}
$$ S_{n}(x):=\sum_{k=1}^{n}b_{k}\sin kx,$$
%\end{equation*}%
and choose $t_{j}=\pi /(2n_{j})$, we have for $k=1,2,\ldots ,n_{j}-1$ that
\begin{eqnarray*}
\sum_{m=4kn_{j}}^{(4k+2)n_{j}-1}b_{m}\sin mt_{j} &\geq &\sum_{n_{j}/2\leq
m\leq 3n_{j}/2}b_{4kn_{j}+m}\sin (4kn_{j}+m)\frac{\pi }{2n_{j}} \\
&\geq &\frac{\sqrt{2}}{2}\sum_{n_{j}/2\leq m\leq 3n_{j}/2}b_{4kn_{j}+m} \\
&\geq &\frac{\sqrt{2}}{2}\frac{1}{\sqrt{\log n_{j}}}\frac{n_{j}}{%
(4k+3/2)n_{j}} \\
&=&\sqrt{2}\frac{1}{\sqrt{\log n_{j}}}\frac{1}{8k+3}.
\end{eqnarray*}%
 On the other hand,
%\begin{equation*}
$$ \left| \sum_{m=(4k+2)n_{j}}^{4(k+1)n_{j}-1}b_{m}\sin mx\right|
\leq
\sum_{m=(4k+2)n_{j}}^{4(k+1)n_{j}-1}b_{m}\leq \frac{1}{8\sqrt{\log n_{j}}}%
\frac{2n_{j}}{(4k+2)n_{j}}\leq \frac{1}{2\sqrt{\log
n_{j}}}\frac{1}{8k+4}.$$
%\end{equation*}%
Therefore
\begin{eqnarray*}
S_{n_{j+1}}(t_{j})-S_{n_{j}}(t_{j}) &\geq &\sum_{k=1}^{n_{j}-1}\left(
\sum_{m=4kn_{j}}^{(4k+2)n_{j}-1}b_{m}\sin mt_{j}-\left|
\sum_{m=(4k+2)n_{j}}^{4(k+1)n_{j}-1}b_{m}\sin mx\right| \right) \\
&\geq &\frac{1}{2}\frac{1}{\sqrt{\log n_{j}}}\sum_{k=1}^{n_{j}-1}\frac{1}{%
8k+4} \\
&\geq &\frac{C}{\sqrt{\log n_{j}}}\log n_{j} \\
&\geq &C\sqrt{\log n_{j}},
\end{eqnarray*}%
and the sine series is not uniformly convergent accordingly since $\sqrt{%
\log n_{j}}\rightarrow \infty $ as $j\rightarrow \infty $.
%\end{pf}

\ \

%\begin{pf}
{\bf PROOF.} [Proof of Proposition 3] If $\{a_{n}\}\in \mbox{\rm
GBVS}$, then for any sufficiently large $n$, we have
\begin{eqnarray*}
\sum_{k=n}^{2n}|\Delta a_{k}| &:&=\sum_{k=n}^{2n}|a_{k}-a_{k+1}|\leq C(%
\mbox{\bf A})\max_{n\leq k<n+N_{0}}a_{k}=:C(\mbox{\bf A})a_{k_{n}}, \\
&&%
\begin{array}{ccc}
\;\;\;\;\;\;\;\;\;\;\;\;\;\;\;\;\;\; & \;\;\;\;\;\;\;\;\;\;\;\;\;\;\;\;\;\;%
\;\;\;\;\;\; & n\leq k_{n}<n+N_{0}\leq 2n.%
\end{array}%
\end{eqnarray*}%
For any $k_{n}/2\leq j\leq k_{n}$,
%\begin{equation*}
$$a_{k_{n}}=\sum_{i=k_{n}}^{2j}\Delta a_{i}+a_{2j+1}\leq
\sum_{i=j}^{2j}|\Delta a_{i}|+a_{2j+1}\leq C(\mbox{\bf
A})a_{k_{j}}+a_{2j+1},$$

%\end{equation*}%
thus
%\begin{equation*}
$$a_{k_{n}}\leq \frac{C\left( \mbox{\bf A}\right)
}{k_{n}}\sum_{k_{n}/2\leq j\leq k_{n}}(a_{k_{j}}+a_{2j+1}).$$
%\end{equation*}%
Since $n\leq k_{n}<n+N_{0}$, $j\leq k_{j}<j+N_{0}$, and $k_{j}$ can repeat
at most $N_{0}$ times for $k_{n}/2\leq j\leq k_{n}$, with the above
estimate, we have

%\begin{equation*}
$$a_{k_{n}}\leq \frac{C\left( \mbox{\bf A}\right) N_{0}}{n}%
\sum_{j=[n/2]}^{4n+1}a_{j},$$

%\end{equation*}

which gives us a $\lambda =5$ such that

%\begin{equation*}
$$\sum_{k=n}^{2n}|\Delta a_{k}|\leq C\left( \mbox{\bf A}\right)
a_{k_{n}}\leq \frac{C\left( \mbox{\bf A}\right)
}{n}\sum_{j=[\lambda ^{-1}n]}^{\lambda n}a_{j}.$$

%\end{equation*}%

On the other hand, let $n_{k}=2^{k}$, $k=0,1,\cdots $, set

%\begin{equation*}
$$a_{n}=\left\{
\begin{array}{ll}
0, & n_{k}\leq n<n_{k}+k, \\
b_{n}, & n_{k}+k\leq n<n_{k+1}-k, \\
0, & n_{k+1}-k\leq n<n_{k+1},%
\end{array}%
\right.
$$
%\end{equation*}%
where $b_{1},\cdots ,b_{n_{k}-k},b_{n_{k}+k},b_{n_{k}+k+1},\cdots
,b_{n_{k+1}-k-1},b_{n_{k+1}+k+1},\cdots $ is any decreasing
nonnegative sequence, we can easily check that $\{a_{n}\}$ does
not satisfy (2) (not in NBVS either), but it certainly belongs to
MVBVS.\\
%\end{pf}

The proof of Proposition 4 is similar to the proof of Proposition 3. Now we
divide the proof of Theorem 5 into the following three lemmas. For a given
series

%\begin{equation*}
$$\sum_{k=1}^{\infty }a_{k}\sin kx=\lim_{n\rightarrow \infty
}\sum_{k=1}^{n}a_{k}\sin kx,$$

%\end{equation*}
we write

%\begin{equation*}
$$ f\left( x\right) =\sum_{k=1}^{\infty }a_{k}\sin kx$$
%\end{equation*}
for those points $x$ where the series converges and let $S_{n}\left(
f,x\right) $ be the $n$th partial sum of $f$ at $x$. As the sequence $%
\{a_{n}\}$ under consideration in the sine series starts with
$a_{1},$ we assume, without loss of generality, that $a_{0}=0.$\\

%\begin{lem}
{\bf Lemma 7}\quad  Let $\{a_{n}\}$ be a nonnegative sequence and
let $f\left( x\right) =\sum_{n=1}^{\infty }a_{n}\sin nx\in C_{2\pi
}.$ Then
%\begin{equation*}
$$\lim_{n\rightarrow \infty }\left| \left| f-S_{n}\left( f\right)
\right| \right| =0.$$
%\end{equation*}
%\end{lem}

This lemma is a direct corollary of Lemma 3 in [7].\\

%\begin{lem}
{\bf Lemma 8}\quad Let $\{a_{n}\}\in \mbox{\rm MVBVS}$. Then
either the uniform convergence of
series (\ref{a1}), or the continuity of its sum function $f$, implies that $%
\lim\limits_{n\rightarrow \infty }na_{n}=0$.\\
%\end{lem}

%\begin{pf}
{\bf PROOF.} Either the uniform convergence of series (\ref{a1}),
or the continuity of its sum function, implies that (by Lemma 7)
%\begin{equation*}
$$\lim_{n\rightarrow \infty }\Vert S_{[\lambda n]}\left( f\right) -S_{[\frac{n%
}{2\lambda }]-1}\left( f\right) \Vert =0$$
%\end{equation*}
holds for any given $\lambda \geq 2$. Since $\{a_{n}\}\in \mbox{\rm MVBVS}$,
there exists a $\lambda \geq 2$ such that for any integer $n>0,$
%\begin{equation*}
$$\sum_{k=n}^{2n}|\Delta a_{k}|\leq \frac{C(\mbox{\bf
A})}{n}\sum_{k=[\lambda ^{-1}n]}^{[\lambda n]}a_{k}.$$
%\end{equation*}
So for $j=n+1,\ldots 2n,$ we have%
\begin{equation}
\label{c3}
\begin{array}{lll}
a_{n} &\leq &\sum\limits_{k=n}^{j-1}\left| \Delta a_{k}\right| +a_{j}\\
&\leq &\sum\limits_{k=\left[ \frac{j}{2}\right] }^{2\left[
\frac{j}{2}\right]
}\left| \Delta a_{k}\right| +a_{j}   \\
&\leq &\frac{C\left( \mbox{\bf A}\right) }{n}\sum\limits_{k=\left[ \frac{j}{%
2\lambda }\right] }^{\left[ \lambda j/2\right] }a_{k}+a_{j} \\
&\leq &\frac{C\left( \mbox{\bf A}\right) }{n}\sum\limits_{k=\left[ \frac{n}{%
2\lambda }\right] }^{\left[ \lambda n\right] }a_{k}+a_{j}.
\end{array}
\end{equation}

Taking the sum of the $n$ inequalities of (\ref{c3}) for $j$ runs
from $n+1$ to $2n,$ we have
\begin{equation}
na_{n}\leq C\left( \mbox{\bf A}\right) \sum_{k=\left[ \frac{n}{2\lambda }%
\right] }^{\left[ \lambda n\right] }a_{k}+\sum_{j=n+1}^{2n}a_{j}\leq C\left( %
\mbox{\bf A}\right) \sum_{k=\left[ \frac{n}{2\lambda }\right] }^{\left[
\lambda n\right] }a_{k}.  \label{c4}
\end{equation}%
Please note that $C\left( \mbox{\bf A}\right) $ may have different values in
different occurrences. Now let $t_{n}=\pi /(2\lambda n)$. Then we have
%\begin{equation*}
$$S_{[\lambda n]}(f,t_{n})-S_{[\frac{n}{2\lambda }]-1}(f,t_{n})\geq C\left( %
\mbox{\bf A}\right) \sum_{k=[\frac{n}{2\lambda }]}^{[\lambda
n]}a_{k}\geq C\left( \mbox{\bf A}\right) na_{n}$$
%\end{equation*}
and the required result follows.\ \ \ \ \
%\end{pf}

%\begin{lem}

{\bf Lemma 9} \quad {\it Let $\{a_{n}\}\in \mbox{\rm MVBVS}$. Then
$\lim\limits_{n\rightarrow \infty }na_{n}=0$ implies that 
$\lim\limits_{n\rightarrow \infty }\Vert f-S_{n}\left( f\right)
\Vert =0$.}\\
%\end{lem}

%\begin{pf}
{\bf Proof.}   We need only to show that
\begin{equation}
\lim_{n\rightarrow \infty }\Vert I(x)\Vert :=\lim_{n\rightarrow \infty
}\left\| \sum_{k=n}^{\infty }a_{k}\sin kx\right\| =0.  \label{c5}
\end{equation}%
In view of $I(0)=I(\pi )=0$, we may restrict $x$ within $(0,\pi )$. From the
condition, for any given $\varepsilon >0$, there is a $n_{0}>0$ such for all
$n\geq n_{0}$ that $na_{n}<\varepsilon $. Let $n\geq \lambda n_{0}$, where $%
\lambda \geq 2$ is the number obtained from the given sequence $\{a_{n}\}\in %
\mbox{\rm MVBVS},$ by the definition for MVBVS. Take $N=[1/x]$ and set
%\begin{equation*}
$$I(x)=\sum_{k=n}^{N-1}a_{k}\sin kx+\sum_{k=N}^{\infty }a_{k}\sin
kx=:I_{1}(x)+I_{2}(x),$$
%\end{equation*}
where, without loss of generality, we assume that $N>n,$ if $N\leq n$, the
same argument as in estimating $I_{2}$ can be applied to deal with $%
\sum_{k=n}^{\infty }a_{k}\sin kx$ directly.\textit{\textrm{\textbf{\ }}}%
Obviously,
%\begin{equation*}
$$|I_{1}(x)|\leq x\sum_{k=n}^{N-1}ka_{k}\leq x(N-1)\varepsilon \leq
\varepsilon .$$
%\end{equation*}
By the well-known estimate
%\begin{equation*}
$$|D_{n}(x)|=\left| \sum_{k=1}^{n}\sin kx\right| \leq \frac{\pi
}{x},$$
%\end{equation*}
and by Abel's transformation and that $\{a_{n}\}\in \mbox{\rm MVBVS}$,
\begin{eqnarray*}
|I_{2}(x)| &=&\left| \sum_{k=N}^{\infty }a_{k}\sin kx\right| \\
&\leq &\sum_{k=N}^{\infty }|\Delta a_{k}||D_{k}(x)|+a_{N}|D_{N-1}(x)| \\
&\leq &Cx^{-1}\left( \sum_{k=N}^{\infty }|\Delta a_{k}|+a_{N}\right) \\
&\leq &CN\sum_{k=N}^{\infty }|\Delta a_{k}|+\varepsilon .
\end{eqnarray*}%
We check that
\begin{eqnarray*}
\sum_{k=N}^{\infty }|\Delta a_{k}| &=&\sum_{j=0}^{\infty }\sum_{2^{j}N\leq
k<2^{j+1}N}|\Delta a_{k}| \\
&\leq &C\left( \mbox{\bf A}\right) \sum_{j=0}^{\infty }\frac{1}{2^{j}N}%
\sum_{k=[\lambda ^{-1}2^{j}N]}^{[\lambda 2^{j}N]}a_{k} \\
&\leq &C\left( \mbox{\bf A}\right) N^{-1}\varepsilon \sum_{j=0}^{\infty
}2^{-j}\sum_{k=[\lambda ^{-1}2^{j}N]}^{[\lambda 2^{j}N]}k^{-1} \\
&\leq &C\left( \mbox{\bf A},\lambda \right) N^{-1}\varepsilon
\sum_{j=0}^{\infty }2^{-j} \\
&\leq &C\left( \mbox{\bf A},\lambda \right) N^{-1}\varepsilon ,
\end{eqnarray*}%
since $\{a_{n}\}\in \mbox{\rm MVBVS}$. Combining all the above
estimates, we have the required result.\\
%\end{pf}

%\begin{pf}
{PROOF.} [Proof of Theorem 6] The construction is to be processed
in a similar but more delicate way to
the proof of Theorem 1. Without loss of generality, we can assume that $%
M_{1}\geq 10$, therefore $M_{j}\geq 10$ for all $j\geq 1.$ Set $n_{1}=1$, $%
n_{2}=10$, and $n_{j+1}=2[M_{4n_{j}}^{1/2}]n_{j}$ for $j=2,3,....$ Let
%\begin{equation*}
$$a_{k}=1,\;\;1\leq k<40.$$
%\end{equation*}
For $j\geq 2$ and $k=1,2,\cdots ,2[M_{4n_{j}}^{1/2}]-1,$ let%
\begin{eqnarray*}
a_{m} &=&\frac{1}{\sqrt{\log M_{4n_{j}}}}\frac{1}{m},\;\;\;if
4kn_{j}\leq m<(4k+2)n_{j}, \\
a_{m} &=&\frac{1}{8\sqrt{\log M_{4n_{j}}}}\frac{1}{m},\;\;if
(4k+2)n_{j}\leq m<4(k+1)n_{j}.
\end{eqnarray*}
Define accordingly a sine series $\sum\limits_{m=1}^{\infty
}a_{m}\sin mx,$ we will show that this series is exact what
required to prove Theorem 6. For any given $n$, there exists a
$j\geq 2$ and a $k,\;1\leq k\leq
2[M_{4n_{j}}^{1/2}]-1$, such that $4kn_{j}\leq n<4(k+1)n_{j}$, then $%
8kn_{j}\leq 2n<8(k+1)n_{j}$. Divide the argument into two cases.\newline
\ \newline
Case 1. $1\leq k\leq \lbrack M_{4n_{j}}^{1/2}]-1$. Then
%\begin{equation*}
$$2n\leq 8[M_{4n_{j}}^{1/2}]n_{j}=4n_{j+1}.$$
%\end{equation*}

We check that
\begin{equation}
\label{c6}
\begin{array}{lll}
\sum\limits_{m=n}^{2n}|\Delta a_{m}| &\leq
&C\sum\limits_{m=k}^{2(k+1)}a_{4mn_{j}}\\
&\leq &\frac{C}{\sqrt{\log
M_{4n_{j}}}}\sum\limits_{m=k}^{2(k+1)}\frac{1}{4mn_{j}}\\
&\leq &\frac{C}{\sqrt{\log M_{4n_{j}}}}\frac{k+3}{4kn_{j}} \\
&\leq &\frac{C}{\sqrt{\log M_{4n_{j}}}}\frac{1}{n_{j}}.
\end{array}
\end{equation}
At the same time,
\begin{equation}
\label{c7}
\begin{array}{lll}
\sum\limits_{m=n}^{2n}a_{m} &\geq &\frac{1}{8\sqrt{\log M_{4n_{j}}}}%
\sum\limits_{m=4(k+1)n_{j}}^{8kn_{j}}\frac{1}{m}   \\
&\geq &\frac{1}{8\sqrt{\log
M_{4n_{j}}}}\frac{(4k-3)n_{j}}{8kn_{j}} \\
&\geq &\frac{1}{64\sqrt{\log M_{4n_{j}}}}.
\end{array}%
\end{equation}
Thus, by noting that $4n_{j}\leq 4kn_{j}\leq n\leq \left( 4k+1\right) n_{j}$%
, $k\leq \lbrack M_{4n_{j}}^{1/2}]-1,$ for any $\lambda \geq 2$, with (\ref%
{c6}) and (\ref{c7}), we have
%\begin{equation*}
$$\frac{\sum_{m=n}^{2n}|\Delta
a_{m}|}{\frac{M_{n}}{n}\sum_{m=[\lambda
^{-1}n]}^{[\lambda n]}a_{m}}\leq \frac{\sum_{m=n}^{2n}|\Delta a_{m}|}{\frac{%
M_{n}}{n}\sum_{m=n}^{2n}a_{m}}\leq C\frac{kn_{j}}{M_{n}n_{j}}\leq C\frac{%
M_{4n_{j}}^{1/2}}{M_{n}}\leq CM_{n}^{-1/2},$$
%\end{equation*}%
and the last quantity in the above inequalities obviously tends to zero as $%
n\rightarrow \infty $.\newline
\ \newline
Case 2. $[M_{4n_{j}}^{1/2}]\leq k<2[M_{4n_{j}}^{1/2}]$. Similarly, we
calculate for this case that (note that $2n\leq
16[M_{4n_{j}}^{1/2}]n_{j}<8n_{j+1}$)
\begin{eqnarray*}
\sum_{m=n}^{2n}|\Delta a_{m}| &\leq &\frac{C}{\sqrt{\log M_{4n_{j}}}}%
\sum_{m=k}^{2[M_{4n_{j}}^{1/2}]-1}\frac{1}{4mn_{j}}+\frac{C}{\sqrt{\log
M_{4n_{j+1}}}}\frac{1}{4n_{j+1}} \\
&\leq &\frac{C}{\sqrt{\log M_{4n_{j}}}}n_{j}^{-1}.
\end{eqnarray*}%
On the other hand, by noting that $[n/2]\leq 2(k+1)n_{j}\leq
4[M_{4n_{j}}^{1/2}]n_{j}$ we achieve that
%\begin{equation*}
$$\sum_{m=[n/2]}^{2n}a_{m}\geq \frac{1}{8\sqrt{\log M_{4n_{j}}}}%
\sum_{m=4[M_{4n_{j}}^{1/2}]n_{j}}^{8[M_{4n_{j}}^{1/2}]n_{j}-1}\frac{1}{m}%
\geq \frac{C}{\sqrt{\log M_{4n_{j}}}}.$$
%\end{equation*}%
Therefore, for any $\lambda \geq 2$, it follows that
%\begin{equation*}
$$\frac{\sum_{m=n}^{2n}|\Delta
a_{m}|}{\frac{M_{n}}{n}\sum_{m=[\lambda
^{-1}n]}^{[\lambda n]}a_{m}}\leq \frac{\sum_{m=n}^{2n}|\Delta a_{m}|}{\frac{%
M_{n}}{n}\sum_{m=\left[ n/2\right] }^{2n}a_{m}}\leq C\frac{%
[M_{4n_{j}}^{1/2}]n_{j}}{M_{n}n_{j}}\leq
C\frac{M_{4n_{j}}^{1/2}}{M_{n}}\leq CM_{n}^{-1/2}.$$
%\end{equation*}%
\newline
\newline
Combining these two cases, in any circumstance, for given $\lambda \geq 2$
we have proved
%\begin{equation*}
$$\lim_{n\rightarrow \infty }\frac{\sum_{k=n}^{2n}|\Delta a_{k}|}{\frac{M_{n}}{%
n}\sum_{k=[\lambda ^{-1}n]}^{[\lambda n]}a_{k}}=0.$$
%\end{equation*}%
In a similar argument to the proof of Theorem 1, by taking $t_{j}=\pi
/(2n_{j})$, we can prove that
%\begin{equation*}
$$S_{n_{j+1}}(f,t_{j})-S_{n_{j}}(f,t_{j})\geq \frac{C}{\sqrt{\log M_{4n_{j}}}}%
\sum_{k=1}^{[M_{4n_{j}}^{1/2}]}\frac{1}{8k+4}\geq C\sqrt{\log
M_{4n_{j}}},$$
%\end{equation*}%
with an observation that $na_{n}\rightarrow 0$ as $n\rightarrow
\infty $, which closely depends on $M_{n}\rightarrow \infty $ as
$n\rightarrow \infty . $ \ So we conclude that the series
constructed is not uniformly convergent although
$na_{n}\rightarrow 0$ as $n\rightarrow \infty $.
%\end{pf}

\section{Results in Complex Space}

Given a trigonometric series $\sum_{k=-\infty }^{\infty
}c_{k}e^{ikx}:=\lim\limits_{n\rightarrow \infty }\sum_{k=-n}^{n}c_{k}e^{ikx}$%
, write
%\begin{equation*}
$$f(x)=\sum_{k=-\infty }^{\infty }c_{k}e^{ikx}$$
%\end{equation*}%
for those points $x$ where the series converges. Denote its $n$th partial
sum $\sum_{k=-n}^{n}c_{k}e^{ikx}$ again by $S_{n}(f,x)$. Define the set%
%begin{equation*}
$$K(\theta _{0}):=\{z:|\arg z|\leq \theta _{0},\;\theta _{0}\in
\lbrack 0,\pi /2)\}.$$
%\end{equation*}

\

It is of great interest to establish an analogue to the Chaundy
and Jolliffe theorem in complex spaces since this will include
sine and cosine series as two particular cases. Previous results
concerning the generalization of Chaundy-Jollif Theorem to complex
space can be found in [7], [23], and [25], etc. In this section we
establish
the following\\

%\begin{thm}
{\bf Theorem 10}\quad {\it Let $\mbox{\bf C}=\{c_{n}\}$ be a
complex sequence satisfying
\begin{equation}
c_{n}\in K(\theta _{0}) and c_{n}+c_{-n}\in K(\theta
_{0}),\;n=1,2,... \label{d1}
\end{equation}%
for some $\theta _{0}\in \lbrack 0,\pi /2)$. If there is a $\lambda \geq 2$
such that
\begin{equation}
\sum_{k=n}^{2n}|c_{k}-c_{k+1}|\leq \frac{C\left( \mbox{\bf C}\right) }{n}%
\sum_{k=[\lambda ^{-1}n]}^{[\lambda n]}|c_{k}|  \label{d2}
\end{equation}%
holds for all $n=1,2,\cdots $, then the necessary and sufficient conditions
for $f\in C_{2\pi }$ and $\lim\limits_{n\rightarrow \infty }\Vert
f-S_{n}\left( f\right) \Vert =0$ are that
\begin{equation}
\lim_{n\rightarrow \infty }nc_{n}=0  \label{d3}
\end{equation}%
and
\begin{equation}
\sum_{n=1}^{\infty }|c_{n}+c_{-n}|<\infty .  \label{d4}
\end{equation}
}
%\end{thm}

Note that the condition (\ref{d1}) in Theorem 10 is weaker than
the analogue condition (6) of Theorem 1 in Le and Zhou [7]. The
proof of Theorem 10 is the result of following
four lemmas.\\

{\bf Lemma 11}\quad{\it
(Xie and Zhou [23], Lemma 2). If a complex sequence $\mbox{\bf C}%
=\{c_{n}\}$ satisfies (\ref{d1}) for some $\theta _{0}\in \lbrack 0,\pi /2)$%
, then $f\in C_{2\pi }$ implies (\ref{d4}).}\\
%\end{lem}

{\bf Lemma 12} \quad{\it If a complex sequence $\{c_{n}\}$
satisfies (\ref{d1}) for some $\theta _{0}\in \lbrack 0,\pi /2)$,
then there is a constant $C\left( \theta _{0}\right) >0$ depending
only on $\theta _{0},$ such that
%\begin{equation*}
$$Rec_{n}\leq \left| c_{n}\right| \leq C\left( \theta _{0}\right) Rec_{n},\;\;n=0,1,2,\ldots
.$$}
%\end{equation*}
%\end{lem}

The proof of this lemma is very straightforward.\\

{\bf Lemma 13}\quad{\it If a complex sequence $\{c_{n}\}$
satisfies (\ref{d1}) for some $\theta _{0}\in \lbrack 0,\pi /2)$
and (\ref{d2}), then $\lim\limits_{n\rightarrow \infty }\Vert
f-S_{n}\left( f\right) \Vert =0$ implies (\ref{d3}).}\\
%\end{lem}

{\bf PROOF.}  As
\begin{eqnarray*}
S_{[\lambda n]}(f,x)-S_{[\frac{n}{2\lambda }]-1}(f,x) &=&\sum_{k=[\frac{n}{%
2\lambda }]}^{[\lambda n]}\left( c_{k}e^{ikx}+c_{-k}e^{-kx}\right) \\
&=&\sum_{k=[\frac{n}{2\lambda }]}^{[\lambda n]}c_{k}\left(
e^{ikx}-e^{-ikx}\right) +\sum_{k=[\frac{n}{2\lambda }]}^{[\lambda n]}\left(
c_{k}+c_{-k}\right) e^{-ikx},
\end{eqnarray*}%
so
\begin{equation}
\left| \sum_{k=[\frac{n}{2\lambda }]}^{[\lambda n]}c_{k}\left(
e^{ikx}-e^{-ikx}\right) \right| \leq \left\| S_{[\lambda n]}\left( f\right)
-S_{[\frac{n}{2\lambda }]-1}\left( f\right) \right\| +\sum_{k=[\frac{n}{%
2\lambda }]}^{[\lambda n]}\left| c_{k}+c_{-k}\right| .  \label{d5}
\end{equation}%
On the other hand, if we let $x=x_{0}=\pi /(2\lambda n),$%
\begin{equation}
\label{d6}
\begin{array}{lll}
\left| \sum\limits_{k=[\frac{n}{2\lambda }]}^{[\lambda
n]}c_{k}\left(
e^{ikx_{0}}-e^{-ikx_{0}}\right) \right| &\geq &\left| \sum\limits_{k=[\frac{n}{%
2\lambda }]}^{[\lambda n]}Rec_{k}\left(
e^{ikx_{0}}-e^{-ikx_{0}}\right) \right|  \\
&=&2\sum\limits_{k=[\frac{n}{2\lambda }]}^{[\lambda n]}Rec_{k}\sin
kx_{0}\\
 &\geq &C\sum\limits_{k=[\frac{n}{2\lambda }]}^{[\lambda
n]}Rec_{k}.
\end{array}
\end{equation}
Now by Lemma 12 and a similar calculation to (\ref{c3}) in the proof of
Lemma 8, for $j=n+1,\ldots ,2n,$
%\begin{equation*}
$$\left| c_{n}\right| \leq \frac{C\left( \mbox{\bf C}\right) }{n}\sum_{k=[%
\frac{n}{2\lambda }]}^{[\lambda n]}\left| c_{k}\right| +\left| c_{j}\right|
\leq C\left( \mbox{\bf C},\lambda ,\theta _{0}\right) \left( \frac{1}{n}%
\sum_{k=[\frac{n}{2\lambda }]}^{[\lambda n]}Rec_{k}+Re%
c_{j}\right) ,$$
%\end{equation*}%
and with (\ref{d6}), (\ref{d5}), and a similar calculation to
(\ref{c4}) in the proof of Lemma 8,
\begin{eqnarray*}
n\left| c_{n}\right| &\leq &C\left( \mbox{\bf C},\lambda ,\theta _{0}\right)
\left( \sum_{k=[\frac{n}{2\lambda }]}^{[\lambda n]}Re%
c_{k}+\sum_{j=n+1}^{2n}Rec_{j}\right) \\
&\leq &C\left( \mbox{\bf C},\lambda ,\theta _{0}\right) \left| \sum_{k=[%
\frac{n}{2\lambda }]}^{[\lambda n]}c_{k}\left(
e^{ikx_{0}}-e^{-ikx_{0}}\right) \right| \\
&\leq &C\left( \mbox{\bf C},\lambda ,\theta _{0}\right) \left( \left\|
S_{[\lambda n]}\left( f\right) -S_{[\frac{n}{2\lambda }]-1}\left( f\right)
\right\| +\sum_{k=[\frac{n}{2\lambda }]}^{[\lambda n]}\left|
c_{k}+c_{-k}\right| \right) ,
\end{eqnarray*}%
then this and Lemma 11 imply (\ref{d3}) if
$\lim\limits_{n\rightarrow \infty }\Vert f-S_{n}\left( f\right)
\Vert =0.$\\
%\end{pf}

%\begin{lem}
{\bf Lemma 14}\quad {\it If a complex sequence $\{c_{n}\}$
satisfies the conditions (\ref{d1}) (for some $\theta _{0}\in
\lbrack 0,\pi /2)$) to (\ref{d4}), then
%\begin{equation*}
$$\lim\limits_{n\rightarrow \infty }\Vert f-S_{n}\left( f\right)
\Vert =0.$$}
%\end{equation*}
%\end{lem}

{\bf PROOF.}  Similar to the proof of the following identity
(\ref{d7}) under the conditions of (\ref{d3}) and (\ref{d4}), we
can easily see that the sequence $\left\{ S_{n}\left( f,x\right)
\right\} $ is a Cauchy sequence for each $x$
and therefore it converges at each $x.$ So we only need to show that%
\begin{equation}
\lim\limits_{n\rightarrow \infty }\left\| \sum_{k=n}^{\infty }\left(
c_{k}e^{ikx}+c_{-k}e^{-ikx}\right) \right\| =0.  \label{d7}
\end{equation}%
For any given $\varepsilon >0,$ from the conditions of (\ref{d3}) and (\ref%
{d4}), there exists a $n_{0}>0,$ such that for all $n\geq n_{0},$ we have%
\begin{equation}
n\left| c_{n}\right| <\varepsilon  \label{d8}
\end{equation}%
and%
\begin{equation}
\sum_{k=n}^{\infty }\left| c_{k}+c_{-k}\right| <\varepsilon .  \label{d9}
\end{equation}%
For $n\geq n_{0},$ write%
\begin{eqnarray*}
\sum_{k=n}^{\infty }\left( c_{k}e^{ikx}+c_{-k}e^{-ikx}\right)
&=&\sum_{k=n}^{\infty }\left( c_{k}+c_{-k}\right)
e^{ikx}+2i\sum_{k=n}^{\infty }c_{k}\sin kx \\
&=&:I_{1}\left( x\right) +2iI_{2}\left( x\right) .
\end{eqnarray*}%
From (\ref{d9}), we have
%\begin{equation*}
$$\left| I_{1}\left( x\right) \right| <\varepsilon .$$
%\end{equation*}%
Follow the same steps in the proof of (\ref{c5}) for Lemma 9, and using (\ref%
{d8}) and Lemma 12 instead, we have%
%\begin{equation*}
$$\left| I_{2}\left( x\right) \right| <\varepsilon .$$
%\end{equation*}%
This complete the proof of Lemma 14.
%\end{pf}

\end{document}